\documentclass[10pt,twoside,a4paper,reqno]{amsart}
\usepackage{amssymb, amsfonts, amsthm, latexsym, amsmath, amscd, epsf,graphicx}

\theoremstyle{plain}
\newtheorem{theo}           {Theorem}
\newtheorem{pro}            {Proposition}
\newtheorem{coro}           {Corollary}
\newtheorem{lemm}           {Lemma}
\newtheorem{conj}           {Conjecture}
\newtheorem{exam}        {Example}

\theoremstyle{definition}
\newtheorem{pr}              {Problem}
\newtheorem*{ack}            {Acknowledgements}

\theoremstyle{remark}
\newtheorem{rem}             {Remark}

\newenvironment{theorem}{\begin{theo}}{\end{theo}}
\newenvironment{proposition}{\begin{pro}}{\end{pro}}
\newenvironment{corollary}{\begin{coro}}{\end{coro}}
\newenvironment{lemma}{\begin{lemm}}{\end{lemm}}
\newenvironment{remark}{\begin{rem}}{\end{rem}}

\newenvironment{problem}{\begin{pr}}{\end{pr}}

\newcommand \CC {\mathbb C}
\newcommand \ZZ {\mathbb Z}

\newcommand \de {\delta}

\newcommand {\cP} {\mathbb {CP}}
\newcommand {\rP} {\mathbb {RP}}
\newcommand \C {\mathcal C}
\newcommand \D {\mathcal D}

\newcommand {\Ct} {\mathbf {Cat}}

\newcommand \ga {\gamma}
\newcommand \pa {\partial}

\newcommand \si {\sigma}

\newcommand \cO{\mathcal O}
\newcommand \maxid {\mathfrak m}

\begin{document}
    \bigskip
    \title[First steps towards total reality of meromorphic functions]
    {First steps towards total reality\\ of meromorphic functions}

\author[T.~Ekedahl]{Torsten Ekedahl}
\address{Department of Mathematics, Stockholm University, SE-106 91
Stockholm, Sweden}
\email{teke@math.su.se}
\author[B.~Shapiro]{Boris Shapiro}
\address{Department of Mathematics, Stockholm University, SE-106 91
Stockholm, Sweden}
\email{shapiro@math.su.se}
\author[M.~Shapiro]{Michael Shapiro}
\address{Department of Mathematics, Michigan State University, East Lansing,  MI 48824-1027,
USA}
\email{mshapiro@math.msu.edu}

\thanks{M.S. is partially supported by NSF grants DMS-0401178, by the BSF grant 2002375
and by the Institute of Quantum Science, MSU}


\begin{abstract}
It was earlier conjectured by the second and the third authors that
  any rational curve $\ga:\cP^1\to\cP^n$ such that the
inverse images of all its flattening points lie on the real line
$\rP^1\subset \cP^1$  is real algebraic  up to a linear fractional
transformation of the image $\cP^n$,  see \cite {EG}, \cite {KhS}
and \cite {SS}. (By a {\em flattening} point $p$ on $\ga$ we mean
a point at which the Frenet $n$-frame
$(\ga',\ga'',\ldots,\ga^{(n)})$  is degenerate.) Below we extend
this conjecture to  the case of
 meromorphic functions on real algebraic curves of higher genera and settle it for
 meromorphic functions of degrees $2,3$ and several other cases.
\end{abstract}

\maketitle

\centerline{\small{\it To Victor Vassiliev on the occasion of his
fiftieth birthday}}

\section {Introduction}

The  above mentioned conjecture  on total reality for rational curves was
formulated by the authors in a private communication to F.~Sottile in 1993
 and attracted some attention due to its
close connection with  the problem of total reality in Schubert
calculus. At the present moment it is supported by a large number
of partial results and extensive numerical evidence, see \cite
{EG},\cite {EGSV},\cite {KhS},\cite {So1}-\cite{So6}, and \cite
{Ve}.  At the same time the only case of this conjecture  which is
completely settled is the case $n=1$, i.e. the case of the usual
rational functions. (The authors were recently informed by
A.~Eremenko and A.~Gabrielov that they proved the above conjecture
in the  case of plane rational quintics.) Namely the main result of \cite {EG} is as follows.

\begin{theorem}[Theorem 2 of \cite {EG} ] For any given $(2d-2)$-tuple of distinct
real numbers there exist at least
$\Ct_{d}=\frac{1}{d}{{2d-2}\choose{d-1}}$ real rational functions
(considered up to a real M\"obius transformation of the image
$\cP^1$) with these critical points. \label{th:EG}
    \end{theorem}

   The above theorem  together with the statement of L.~Goldberg, see \cite{Go} claiming
    that for any $(2d-2)$-tuple of distinct complex numbers the number
    of complex rational functions (considered up to a complex M\"obius
transformation of the image $\cP^1$) with these critical points is at
most $\Ct_{d}$ gives the proof in the case  $n=1$.
The main idea of the proof of Theorem~\ref{th:EG} is the
explicit construction of such functions using the notion of {\em
garden} which is the graph on the source $\cP^1$ obtained as the inverse image of $\rP^1$
under a real rational function, comp. \cite {NSV}. The number of topologically different
gardens for generic real rational functions of degree $d$ with all
real, simple and distinct critical points turned out
to coincide with $\Ct_{d}$. The difficult part of the proof of
Theorem~\ref{th:EG} is then to show that for a given topological type of a
garden there exists a real rational function with this garden and having
$(2d-2)$ prescribed  real critical points.
\medskip

The purpose of this short paper is to discuss  a (conjectural)
generalization of  Theorem~\ref{th:EG} to the case of the source curves of higher genera,
i.e. to the case of meromorphic functions.
Existence of real meromorphic functions with all real (and closely located) critical points
on real  curves of positive genus was recently proved by B.~Osserman in \cite {Os2}.

We start with some standard notation.

\medskip
\noindent
{\bf Definition.} A pair $(\C,\sigma)$ consisting of a compact Riemann surface $\C$
and its antiholomorphic involution $\sigma$  is called a
{\em real algebraic} curve.
It is well-known that if $\C$ is a compact Riemann surface of genus
$g$ then for any $\sigma$ the set $\C_{\sigma}$ (if nonempty) consists of
at most $g+1$ disjoint smooth closed non-selfintersecting loops called the {\em ovals} of
$\C_{\sigma}$. The set $\C_{\sigma}\subset \C$ of all fixed points
of $\sigma$ is called the {\em real part} of $(\C,\sigma)$.

If $(\C,\sigma)$ and
$(\D,\tau)$ are real curves (varieties) and ${f}:{\C}\to{\D}$
a holomorphic map, then we shall use the notation $\overline{f}$ for $\tau\circ
f\circ\sigma$ which is another holomorphic map. The map $f$ is {\em real}
if $\overline{f}=f$.

\medskip
\noindent
{\bf Definition.} Following the terminology of real algebraic
geometry we call a real algebraic curve $(\C,\sigma)$ with $\C$
compact of genus $g$ an {\em $M$-curve} if its $\C_{\sigma}$ consists
of exactly $g+1$ ovals.

The main question we discuss  below is as follows.

\medskip
\begin{problem} Given a meromorphic function
$f:(\C,\sigma)\to \cP^1$ such that

\noindent
i) all its critical points and values are distinct;

\noindent
 ii) all its critical points belong to $\C_{\sigma}$

\noindent is it true that that $f$ becomes a real meromorphic
function after a choice of a real structure of $\cP^1$?
\end{problem}


\medskip
\noindent
{\bf Definition} We say that the space of meromorphic functions of
degree $d$ on a genus $g$ real algebraic curve $(\C,\sigma)$
has \emph{the total reality property} if Problem~1 has the
affirmative answer for any meromorphic function from this space which   satisfies the above assumptionis.


\medskip
Notice that Problem 1 has the following modification. Since the number of critical
points/values of a generic degree $d$
meromorphic function from a genus $g$ curve equals $2d-2+2g$ one has that the
dimension of the space of corresponding linear systems equals
$2d-2-g$, i.e. one can arbitrarily assign the position of $2d-2-g$
critical points and find (in general, several)  meromorphic functions with these critical points. For each such function the remaining  $3g$ critical points will be uniquely determined.

\medskip
\noindent
\begin{problem}
Given a meromorphic function $f:(\C,\sigma)\to \cP^1$ of degree
$d$ such that

\noindent
i) all its critical points and values are distinct;

\noindent
 ii)  its $2d-2-g$  critical points belong to $\C_{\sigma}$

\noindent is it true that that $f$ becomes a real meromorphic
function after a choice of a real structure of $\cP^1$?
\end{problem}

\medskip
In the present note we prove the following results.

 \begin{theorem}  The space of  meromorphic  functions of  any degree $d$ which is a prime
 on any real curve $(\C,\sigma)$ of genus $g$ which
 additionally satisfies the inequality:
 $g>\frac{d^2-4d+3}{3}$  has the total reality property.
\label{th:sqrt}
\end{theorem}



%
%
%
%

\begin{corollary}  The  total reality property holds for all meromorphic functions
 of  degrees $2,3$. Moreover, if degree
equals $2$ then in case of positive genus it is sufficient that
just one critical point is real. \label{th:main}
\end{corollary}


Relaxing the requirement that all critical points must be distinct
we obtain the following corollary.

\begin{coro} Any meromorphic function of a prime degree $d$ on a real algebraic surface
$(\C,\sigma)$, such that $g(\C)\ge (d-1)^2$ and all its (not
necessary distinct) critical points belong to $\C_\sigma$ becomes
a real function after appropriate M\"obius transformation of the
image $\cP^1$.
    \label{cr:main}
    \end{coro}

The total reality property for degree $4$ meromorphic functions is
shown to be equivalent to the non-existence of a sextic in $\cP^2$
whose only  singularities are:  $7$ real cusps and two complex
conjugated nodes, such that the line connecting two nodes is
tangent to the sextic in a real smooth point. Unfortunately, the
known results about the realizability of singularities by real
algebraic curves are not strong enough to cover  this remaining
case, comp. \cite{Shu}, \cite {Shu2}.



\medskip
 On the other hand, we show that the answer to Problem 2 is negative. Namely,

\begin{proposition} There exists a  real elliptic  curve  $(\C,\sigma)$ with a nonempty real part
    $\C_{\sigma}$ and a meromorphic function $f:\C\to\cP^1$  of degree $3$
    with $3$ of its $6$ critical points lying on  $\C_{\sigma}$ and which can not be made
     real by a M\"obius transformation of the    image $\cP^1$.
    \label{pr:2}
    \end{proposition}


\begin{rem} Note that (as it was pointed out to the second author by A.~Eremenko)
a further generalization of the Problem~1 to the case  of maps
between two real curves definitely has a negative answer. The
counterexample can  be found already for a map between two
elliptic curves. By the Riemann-Hurwitz formula such a map has no
ramification locus. If reality conjecture holds then any map
$f:E_1\to E_2$ from a real elliptic curve $(E_1,\sigma)$ to an
elliptic $E_2$ is real, i.e., there exists an involution
$\tau:E_2\to E_2$ such that $\bar f=\tau f\sigma=f$. However, it
is evidently false. Namely, let $E_1=E_2=(\CC/{\ZZ+i\ZZ},\sigma)$
where $\sigma$ is the standard complex conjugations. Let
$\Phi:\CC\to\CC$ be the linear map $z\mapsto (2+i)z$, and
$\phi:E_1\to E_2$ be degree 5 map induced by $\Phi$. Existence of
$\tau$ means, in particular, that push forward of $\sigma$ is well
defined. Put $\xi=\frac{i}{2+i}=\frac{1+2i}{5}=\Phi^{-1}(i)$.
$\Phi(\bar\xi)=\frac{4-3i}{5}\not\in \ZZ+i\ZZ$, which shows that
the push forward of $\sigma$ under $\phi$ is not well defined.
Hence, there is no antiholomorphic involution on $E_2$ which makes
$\phi$ into a real map.
\end{rem}

\medskip
The structure of the note is as follows. \S~\ref{sc:pr} contains the
proofs of Theorem~\ref{th:main}, Corollary~\ref{cr:main} and Proposition~\ref{pr:2} while
\S~\ref{sc:rmk} contains a number of remarks and open problems.

\begin{ack} The authors are sincerely grateful to
A.~Gabrielov, A.~Eremenko,  R.~Kulkarni, S.~Natanzon, B.~Osserman,
A.~Vainshtein and, especially F.~Sottile
 for numerous discussions of the topic.  The second and the third authors want to
acknowledge the hospitality of MSRI in Spring 2004 during the
program 'Topological methods in real algebraic geometry' which
gave them a large number of valuable research inputs.
    \end{ack}

\section {Proofs}
\label{sc:pr}

We start our proofs with a characterization of  real meromorphic functions on a real
algebraic curve $(\C,\sigma)$.
\begin{proposition}\label{realification criterion}
If $(\C,\sigma)$ is a proper irreducible real curve and ${f}:{\C}\to{\cP^1}$ a
non-constant holomorphic map (with $\cP^1$ provided with its standard real structure), then
$f$ is real for some real structure on $\cP^1$ precisely when there is a M\"obius
transformation ${\varphi}:{\cP^1}\to{\cP^1}$ such that $\overline{f}=\varphi\circ f$.
\begin{proof}
Any real structure on $\cP^1$ is of the form $\tau\circ\phi$ for a
complex M\"obius transformation $\phi$ and $\tau$ the standard
real structure with $\overline{\phi}\circ\phi=\text{id}$ and
conversely any such $\phi$ gives a real structure. If $f$ is real
for such a structure we have $f=\tau\circ \phi\circ f\circ\sigma$,
i.e., $\overline{f}=\varphi\circ f$ for
$\varphi=\tau\circ\phi^{-1}\circ\tau$. Conversely, if
$\overline{f}=\varphi\circ f$, then $f=
\overline{\overline{f}}=\overline{\varphi\circ
f}=\overline{\varphi}\circ\overline{f}=\overline{\varphi}\circ\varphi\circ
f$ and as $f$ is surjective we get
$\overline{\varphi}\circ\varphi=\text{id}$. That means that
$\phi:=\tau\circ\varphi^{-1}\circ\tau$, then $\phi$ defines a real
structure on $\cP^1$ and by construction $f$ is real for that
structure and the fixed one on $\C$.
\end{proof}
\end{proposition}
We recall that up to  a real isomorphism there are only two real structures on
$\cP^1$, the standard one and the one on an isotropic real quadric in $\cP^2$. The
latter is distinguished from the former by not having any real points.

Assume now that $(\C,\sigma)$ is a proper irreducible real curve and ${f}:{\C}\to{\cP^1}$ a
non-constant  holomorphic map. It defines the holomorphic map
\begin{displaymath}
\C \stackrel{(f,\overline{f})}{\longrightarrow} \cP^1\times\cP^1
\end{displaymath}
and if $\cP^1\times\cP^1$ is given the real structure that takes $(x,y)$ to
$(\tau(y),\tau(x))$, which we shall call the {\em involutive real
structure}, then it is clearly a real map.

\noindent
\begin{proposition}\label{real diagonal} \

\begin{enumerate}
\item \label{it1} The image $\D$ of the curve $\C$ under the map $(f,\overline{f})$ is of type
$(\de,\de)$ for some positive integer $\de$ and if $\pa$ is the degree of the map $\C \to
\D$ we have that $d=\de\pa$, where $d$ is the degree of the original $f$.

\item \label{it2} The function $f$ is real for some real structure
on $\cP^1$ precisely when $\de=1$.

\item \label{it3} Assume that $\C$ is smooth and all the critical points of $f$ are real. Then all the
critical points of  ${g}:{\widetilde \D}\to{\cP^1}$, the composite of the normalization
map $\widetilde \D \to \D$ and the restriction of the projection of $\cP^1\times\cP^1$
has all its critical points real.
\end{enumerate}
\end{proposition}

\begin{proof}
The image under the real map $(f,\overline{f})$ of $\C$ is a real curve so that
$\D$ is a real curve in $\cP^1\times\cP^1$ (with its involutive real
structure). Any such curve is of type $(\de,\de)$ as the real structure permutes the
two degrees. The rest of (\ref{it1}) follows by using the multiplicativity of degrees for the
maps $f\colon \C \to \D \to \cP^1$, where the last map is projection on the first
factor.

As for (\ref{it2}) assume first that $f$ can be made real for some real
structure on $\cP^1$. By Proposition \ref{realification criterion} there is a
M\"obius transformation $\varphi$ such that $\overline{f}=\varphi\circ f$ but that
in turn means that $(f,\overline{f})$ maps $\C$ into the graph of $\varphi$ in
$\cP^1\times\cP^1$ and that graph is hence equal to $\D$ and is thus of type
$(1,1)$. Conversely, assume that $\D$ is of type $(1,1)$. Then it is a graph of
an isomorphism $\varphi$ from $\cP^1$ to $\cP^1$ and by construction
$\overline{f}=\varphi\circ f$ so we conclude by another application of
Proposition \ref{realification criterion}.

Finally, for (\ref{it3}) we have that the map $\C \to \D$ factors as a,
necessarily real, map ${h}:{\C}\to{\widetilde \D}$ and then $f=g\circ h$. If $pt \in
\widetilde \D$ is a critical point, then all points of $h^{-1}(pt)$ are critical for
$f$ and hence by assumption real. As $h$ is real this implies that $pt$ is also
real.
\end{proof}

Part (\ref{it2}) of the above Proposition gives another reformulation
of the  total reality property  for meromorphic functions.

\begin{corollary}\label{cor:reform}  If a degree $d$ function $f:(\C,\si)\to \cP^1$
is real for some real structure on $\cP^1$ then the map $\C
\stackrel{(f,\overline{f})}{\longrightarrow} \D\subset \cP^1\times
\cP^1$ must have degree $d$ as well.
\end{corollary}

\begin{remark} Notice that without the requirement of reality of $f$ the degree of
$\C \stackrel{(f,\overline{f})}{\longrightarrow}  \D$ can be any factor of $d$.
\end{remark}

By a {\em cusp} we  mean a  curve singularity of multiplicity $2$
and whose tangent cone is a double line. It has the local form
$y^2=x^{k}$ for some integer $k\ge 3$ where $k$ is an invariant
which we shall call its {\em type}. A cusp of type $k$ gives a
contribution of $\lceil(k-1)/2\rceil$ to the arithmetic genus of a
curve. A cusp of type $3$ will be called {\em ordinary}.

If $\C$ is a curve and $p_1,\dots,p_k$ are smooth points on it
then we let ${\pi}:{\C}\to{\C(p_1,\dots,p_k)}$ be the finite map
which is a homeomorphism and for which $\cO_{\C(p_1,\dots,p_k)}
\to \pi_*\cO_\C$ is an isomorphism outside of $\{p_1,\dots,p_k\}$
with $\cO_{\C(p_1,\dots,p_k),\pi(p_i)} \to \cO_{\C,p_i}$ having
image the inverse image of $\CC$ in $\cO_{\C,p_i}/\maxid_{p_i}^2$.
In other words, $\C(p_1,\dots,p_k)$ has ordinary cusps at points
$\pi(p_i)$.

Then $\pi$ has the following two (obvious) properties:


\begin{lemma}\label{lm:properties} \

\begin{enumerate}
\item A holomorphic map ${f}:{\C}\to{X}$ which is not an immersion at all the points
$p_1,\dots,p_k$ factors through $\pi$.

\item If $\C$ is proper, then the arithmetic genus of $\C(p_1,\dots,p_k)$ is $k$
plus the arithmetic genus of $\C$.
\end{enumerate}
\end{lemma}

\begin{proposition}\label{real classification}
Assume that $(\C,\si)$ is a smooth and proper real curve and let
${f}:~{\C}\to{\cP^1}$ be a holomorphic map of degree $d$. If there
are $k$ real points $p_1,\dots,p_k$ on $\C$ which are critical
points for $f$ and if $(f,\overline{f})$ gives a map of degree $1$
from $\C$ to its image $\D$ in $\cP^1\times\cP^1$, then
$g(\C)+k\le (d-1)^2$. If $g(\C)+k=(d-1)^2$, then the map $h:\C \to
\D$ factors to give an isomorphism $\C(p_1,\dots,p_k)
\stackrel{\sim}{\rightarrow} D$.
\begin{proof}
As $p_i$ is real it is a critical point also for $\overline{f}$ and hence for
$(f,\overline{f})$. This implies by the first property for $\C \to
\C(p_1,\dots,p_k)$ that the map $\C \to D$ factors as $\C \to \C(p_1,\dots,p_k) \to
\D$ and hence the arithmetic genus of $\C(p_1,\dots,p_k)$, which is $g(\C)+k$ by
the second property of $\C(p_1,\dots,p_k)$, is less than or equal to the
arithmetic genus of $\D$, which by the adjunction formula is equal to
$(d-1)^2$. If we have equality then their genera are equal and hence the map
$\C(p_1,\dots,p_k) \to \D$ is an isomorphism.
\end{proof}
\end{proposition}

Now we are ready to start the proof of Theorem~\ref{th:sqrt}.
It is now a simple corollary of Proposition~\ref{real
classification}. Indeed,
 if we assume that all the critical points of a generic meromorphic
 function $f:(\C,\sigma)\to\cP^1$
are real then $k$ in the above Proposition equals $2d-2+2g(\C)$.
Under the assumption $g(\C)>\frac{d^2-4d+3}{3}$ one gets
$g(\C)+k=2d-2+3g(\C)>(d-1)^2$. Thus, the case when $\C$ maps
birationally to $\D$ is impossible by the above Proposition. Since
$d$ is prime the only other possible case is when the degree of
the map $\C\to \D$ equals $d$ and therefore, the degree of the map
$\D\to \cP^1$ equals $1$ which by (\ref{it2}) of
Proposition~\ref{real diagonal} gives the the total reality
property holds.\qed

\medskip

We start with the proof of Corollary~\ref{th:main} for the case of
degree $2$, i.e. in the hyper-elliptic situation. Then we will
reprove and strengthen the same result in a different way. But the
simple classical argument below has an independent interest.

\begin{proof} Consider a degree $2$ meromorphic function
$f:(\C,\sigma)\to \cP^1$. It defines the
holomorphic involution on $\C$ which we abusing notation denote by
$f$ as well. Since $\C$ is smooth then by Riemann-Hurwitz formula
the involution $f$ has exactly $2g+2$ distinct fixed points on
$\C$ which are the critical points of the function $f$. Take the
map $\tilde f=\sigma\cdot f\cdot\sigma$. Then $\tilde f$ is also a
holomorphic involution on  $\C$. It is known that on a given
hyper-elliptic curve of genus $g>1$ there exist exactly one
holomorphic involution, see \cite {GH}, sect. 2.3.  Therefore,
$f=\tilde f$, or, equivalently $f\cdot\sigma=\sigma\cdot f$.
Therefore, the antiholomorphic involution $\sigma$ on $\C$ induces
the antiholomorphic involution $\bar \sigma$ on the quotient
$\C/f\simeq \cP^1$. Moreover, the real part $\C_{\sigma}$ is
projected by $f$ on the real part $\rP^1\subset \cP^1$ w.r.t the
antiholomorphic involution $\bar \sigma$ of the image $\C/f\simeq
\cP^1$. Thus the images of all the critical points of $f$ are real
as well, i.e. belong to the above $\rP^1$. This implies that our
original curve $\C$ is realized as the Riemann surface of the
plane algebraic curve given by the equation:
$y^2=(x-a_{1})\ldots(x-a_{2g+2})$ with all $a_{i}$'s
 real and the function $f$ coincides with the projection $(y,x)\to x$ where
$(y,x)$ satisfies the above equation. Notice that in this argument
we also used the fact that any antiholomorphic involution on
$\cP^1$ possessing a fixed point is conjugate to the standard
complex conjugation on $\cP^1$, see Proposition~\ref{realification
criterion}.
 The fact that $(\C,\sigma)$ is an $M$-curve follows
from the consideration of the above equation with all real
$a_{i}$'s. In the case $g=1$  we also get the second holomorphic
involution $\tilde f=\sigma\cdot f\cdot\sigma$. Notice that
$\tilde f$ and $f$ have the same set of $4$ fixed points. It is
well-known that for any point $p$ on $\C$ of genus $g=1$ there
exists and unique holomorphic involution on $\C$ having $p$ as its
fixed point. Therefore, $\tilde f=f$ and the rest of the argument
applies.
    \end{proof}

\medskip
Let us now apply Theorem~\ref{th:sqrt} to prove
Corollary~\ref{th:main} in its full generality.

\medskip
\noindent Case $d=2$. Suppose that the degree $d$ of the map
${f}:(\C,\si)\to{\cP^1}$ is equal to $2$. That only leaves two
possibilities: The first is that the map $\C \to \D$ has degree
$2$ and then by Proposition \ref{real diagonal} $f$ is real for
some real structure $(\cP^1,\tau)$ on $\cP^1$. In particular, if
the set $\C_\si$ of real points is nonempty then $(\cP^1,\tau)$
has the same property which means that it is equivalent to the
standard real structure. The second is that the map $\C \to \D$ is
birational and then by Proposition \ref{real classification} we
get $g(\C)+k \le 1^2=1$, where $k$ is the number of real critical
points of $f$. In particular if $g(\C)>0$ then there are no real
critical points. Thus a hyper-elliptic map from a real curve
$(\C,\si)$ is real if one of its critical points is real.

\noindent
Case $d=3$.  In this case again we have only two possibilities; either $f$
is real for a real structure on $\cP^1$ or $\C \to \D$ is birational in which case
we have $g(\C)+k \le 2^2=4$. The case $g(\C)=0$ was settled in \cite {EG}. Recall that
the total number of critical points equals $2d-2+2g(\C)$. But if $g(\C)>0$
then $2\cdot 3-2+3g>4$ and this case of Theorem~\ref{th:main} is settled.
Analogously to the case $d=2$  a function $f$  with the degree  $d=3$ is real if it has more
than $\max(4-g(\C),1)$ real critical points. \qed

\medskip
Next we prove Corollary~\ref{cr:main}.

\begin{proof} Notice that since $d$ is prime again it is enough to
eliminate the case $\deg h=1$. Then as in proof of the main
theorem the arithmetic genus of $\C(p_1,\dots,p_k)$ must not
exceed the arithmetic genus of $\D$. The arithmetic genus of $\D$
equals $(d-1)^2$ whereas the arithmetic genus of
$\C(p_1,\dots,p_k)$ is the sum of $g(\C)$ and a contribution from
the singular points $p_1,\dots,p_k$. The latter is at least $k\ge
\frac{2d+2g(\C)-2}{d-1}\ge 2$. Hence, the hypothesis $g(\C)\ge
(d-1)^2$ contradicts the assumption $\deg h=1$.
\end{proof}

The case $d=4$ can be reduced to the problem of existence for some
special plane sextics. Indeed, we have three possibilities; $\D$
has degree $(1,1)$, $(2,2)$, or $(4,4)$. In the first case $f$ can
be made real. In the second case, by Proposition \ref{real
diagonal}, the projection on the first factor will give a map from
the normalization $\widetilde \D$ of $\D$. The arithmetic genus
$p_a(\D)=1$, and the geometric genus $g(\widetilde \D)$ of the
normalization $\widetilde D$ does not exceed $1$. Let $\widetilde
h:\C\to\widetilde\D$ be the lift of $h:\C\to\D$. Note that if
$p_i\in \C$ is a critical point of $f$ then either its image
$h(p_i)$ is a cusp of $\D$ or $p_i$ is a ramification point of
$\widetilde h$. The ramification divisor $R(\widetilde
h)=2g(\C)+2-4g(\widetilde\D)$. The number of cusps of $\D$ does
not exceed $1$, whereas the number of distinct critical points of
$f$ is $2g(\C)+6$. Note also that any preimage of a cusp of $\D$
must be a critical point. Since $\deg h=2$, we must have
$\frac{1}{2}\left(2g(\C)+6-\left(2g(\C)+2-4g(\widetilde\D)\right)\right)\le
1$ which is impossible.

 We are hence left with the case when $\D$ has degree
$(4,4)$. The only case when $2\cdot 4-2+3g(\C)\le 9$ for $g(\C)>0$
is the case  of $g(\C)=1$. If all the critical points
$p_1,\dots,p_8$ of  ${f}:{\C}\to{\cP^1}$ are real, then we get a
birational map $\C(p_1,\dots,p_8) \to \D$ and as then both
$\C(p_1,\dots,p_8)$ and $\D$ have arithmetic genus $9$, this map
is an isomorphism. Hence $\D$ is a curve with $8$ ordinary real
cusps and no other singularities. Projecting from a cusp $p$ gives
us a  real plane curve $\D'$ of degree $6$. Apart from the $7$
surviving real cusps, $\D'$ will have a pair of complex conjugate
singularities obtained by contracting the strict transforms of the
two lines through $p$. These strict transforms will have
intersection number two with the strict transform of the curve and
hence the singularities will be ordinary nodes (when they
intersect in two points) or cusps (when they are tangent). The
line through these two points will then be a tangent to $\D'$.
Conversely, if we have an irreducible plane curve of degree $6$
with $7$ real ordinary cusps and two complex conjugate ordinary
nodes or cusps such that the line through them is tangent to the
curve we can go backwards and get a curve of degree $(4,4)$ with
the desired properties.

Note also that if one composes the map from the plane blown up at
two points to a quadric with the projection onto one of the
rulings one obtains the projection from one of the two blown up
points. Hence, given a plane curve $\D'$ of degree $6$ as above,
the map $f$ is obtained by composing the normalization map $\C \to
\D'$ with the projection from one of the non-real singularities.

\begin{remark}As a curiosity, if one has a, not necessarily real, irreducible
plane curve $\D'$ of degree $6$ with $7$ real ordinary cusps and
one tangent at a smooth real point that intersects $\D'$ in two
complex conjugate singular points, then the curve is real. Indeed,
if not it is distinct from it is complex conjugate
$\overline{\D'}$ and by the Bezout theorem the total intersection
multiplicity $\D'\cdot\overline{\D'}=36$, whereas the nine
singular points of $\D'$ are also singular points of
$\overline{\D'}$ giving a total contribution of at least $9\cdot
4=36$ and the smooth real point gives an extra contribution of at
least $1$ (at least $2$ actually as its tangent is common to both
curves).
\end{remark}

Consider the existence problem for real algebraic curves in
$\cP^2$ of degree $6$ with the following singularities: 7 real
cusps and either 2 complex conjugate cusps or nodes with
additional property that the line through those is tangent to the
real part at a smooth point. The situation with $7$ real and $2$
complex conjugate cusps is easy to reject, comp. \cite{ISh}.
Namely, any complex degree $6$ genus $1$ curve has at most $9$
cusps. If it is additionally real and has exactly $9$ cusps then
it is dual to a nonsingular real curve of degree $3$ with exactly
$3$ real inflection points, see e.g. \cite {Wa}. But then the
original curve has exactly $3$ (and no more) real cusps.

Therefore in order to settle the total reality conjecture for all
degree $4$ meromorphic functions it remains to prove or disprove
the existence of a sextic in $\cP^2$ with $7$ real cusps and $2$
complex conjugate nodes such that the line connecting two nodes is
tangent to the sextic at a smooth real point.

\medskip
  Finally, let us settle Proposition~\ref{pr:2}.
 Suppose that $d=3$, then again we have only two possibilities;
either $f$ is real for a real structure on $\cP^1$ or $\C \to \D$
is birational in which case we have $g(\C)+k \le 2^2=4$. Assume
further that $g(\C)=1$ and $k=3$ so that $\D \subset
\cP^1\times\cP^1$ is a real curve of type $(3,3)$ with three real
cusps. Given an arbitrary M\"obius transformation $\varphi$, the
holomorphic map $\Phi\colon(x,y) \mapsto
(\varphi(x),\overline{\varphi(y)})$ is a real automorphism of
$\cP^1\times\cP^1$ with involutive real structure. Projection on
the first factor identifies the real points of $\cP^1\times\cP^1$
(still with its involutive structure) with $\cP^1$ and under this
identification $\Phi$ acts on the real points as $\varphi$ acts on
$\cP^1$. This means that we may assume that the three cusps are
$(0,0)$, $(1,1)$ and $(\infty,\infty)$.

Now, the complete linear system of type $(1,1)$ embeds
$\cP^1\times\cP^1$ as a quadric $Q$ in $\cP^3$. This system has a
real structure with respect to the involutive real structure and
realizes $\cP^1\times\cP^1$ as a quadric of signature
$(+1,-1,-1,-1)$. It has the property that there are real points on
it but the two rulings on it are not defined over the reals,
instead through each real point on it there are two complex
conjugate lines on $Q$ passing though it. Projecting $Q$ from a
real point $p$ gives a map from $Q$ with $p$ blown up to the plane
which gives an isomorphism from $Q$ with $p$ blown up and the
strict transform of the two lines passing through it blown down,
taking the exceptional curve to a line through the two blown down
curves. Conversely, given two complex conjugate points $q$ and
$\overline{q}$ in the projective plane, the linear system of
quadrics passing through them gives a map from the plane with $q$
and $\overline{q}$ blown up onto a quadric in $\cP^3$ which
contracts exactly the line through $q$ and $\overline{q}$.

Now, assume that $\D \subset \cP^1\times\cP^1$ is a real curve of
type $(3,3)$ with $3$ real ordinary (i.e., of type $3$) cusps. We
now project from one of the cusps $p$. The strict transform
$\tilde\D$ of $\D$ in the blowing up $Q$ of $\cP^1\times\cP^1$
will then meet the exceptional curve in one real point. The strict
transform of the two lines through $p$, $E_1$ and $E_2$ are
complex conjugate in $Q$ and meet $\tilde D$ transversally in one
point each (which are each other's complex conjugate). Mapping to
the projective plane gives a curve $\D'$ of degree $4$ with two
real ordinary cusps and no other singular points as well as one
smooth real point whose tangent intersections $\D'$ in the point
and two complex conjugate points. Conversely, suppose $\D'$ is a
plane curve of degree $6$ with two real ordinary real cusps, no
other singularities, and a smooth real point whose tangent
intersects the curve in the point and two complex conjugate
points. Blowing up those two complex conjugate points and blowing
down the tangent gives a curve on $\cP^1\times\cP^1$ with three
ordinary real cusps.

We will show that there exist curves in $\cP^1\times\cP^1$ real in
the involutive real structure, having degree $(3,3)$ and $3$ real
ordinary cusps and no other singularities. Above we explained that
this is equivalent to constructing a plane curve $\D'$ of degree
$4$ with two real ordinary real cusps, no other singularities, and
a smooth real point whose tangent intersects the curve in that
point and two complex conjugate points. Indeed, blowing up those
two complex conjugate points and blowing down the tangent gives a
curve on $\cP^1\times\cP^1$ with three ordinary real cusps.
\begin{proposition}\label{counter}
The space of degree $4$ plane curves with two ordinary cusps and no other
singularities is a smooth subvariety of the space of all quartics. In the real
locus of this space, the conditions that the two cusps are both real and that
there is a smooth point whose tangent intersects the curve in the point and two
complex conjugated points is open and nonempty.
\begin{proof}
For the non-real part it will be enough to show that for a given
curve $\D=\{f=0\}$ in the space, the map from the linear space of
quartics to the product of the tangent spaces to mini-versal
deformations of the two singularities is surjective. That product
can be identified with the product of the ``Milnor spaces''
$\cO_{D,p}/(f_x,f_y,f_z)$ at the two singular points. As the
points are ordinary cusps we have that $\maxid_{\cP^2,p}^2\subset
(f_x,f_y,f_z)$ but it is clear that the quartics fill out
$\cO_{\cP^2,p}/\maxid_p^2\oplus\cO_{\cP^2,q}/\maxid_q^2$ for any
two points of $p$ and $q$ of $\cP^2$.

For the real part, the openness is clear and hence it is enough to
show that it is nonempty. It is easily verified that
$f=xz^3-yz^3+xyz^2+x^2y^2$ has ordinary cusps at $(1:0:0)$ and
$(0:1:0)$ and no other singularities. $(0:0:1)$ lies on the curve
and its tangent is given by $x=y$ which intersects the curve in
$(0:0:1)$, $(i:i:1)$, and $(-i:-i:1)$.


%
The latter statement implies that there exist real curves of
bidegree $(3,3)$ with exactly $3$ real ordinary cusps and no other
singularities which settles Proposition~\ref{pr:2}.
\end{proof}
\end{proposition}

A more explicit example of a degree $3$ function from a real elliptic curve with $3$
real critical points but which can not be made real is presented below.


\begin{figure}[!htb]
\centerline{\hbox{\epsfysize=4cm\epsfbox{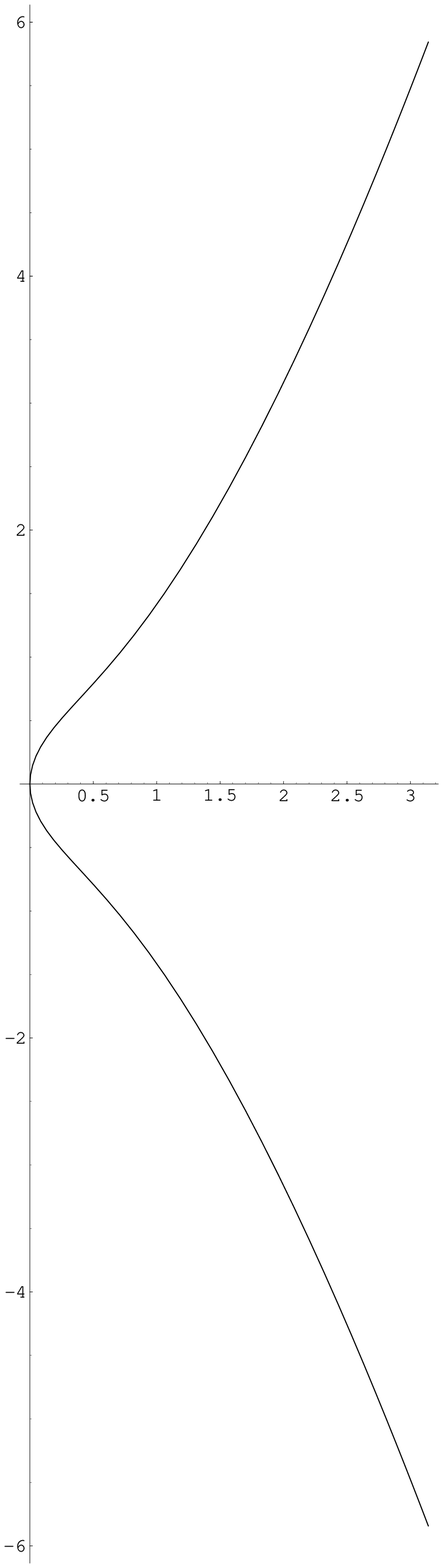}}
\hskip0.5cm\hbox{\epsfysize=4cm\epsfbox{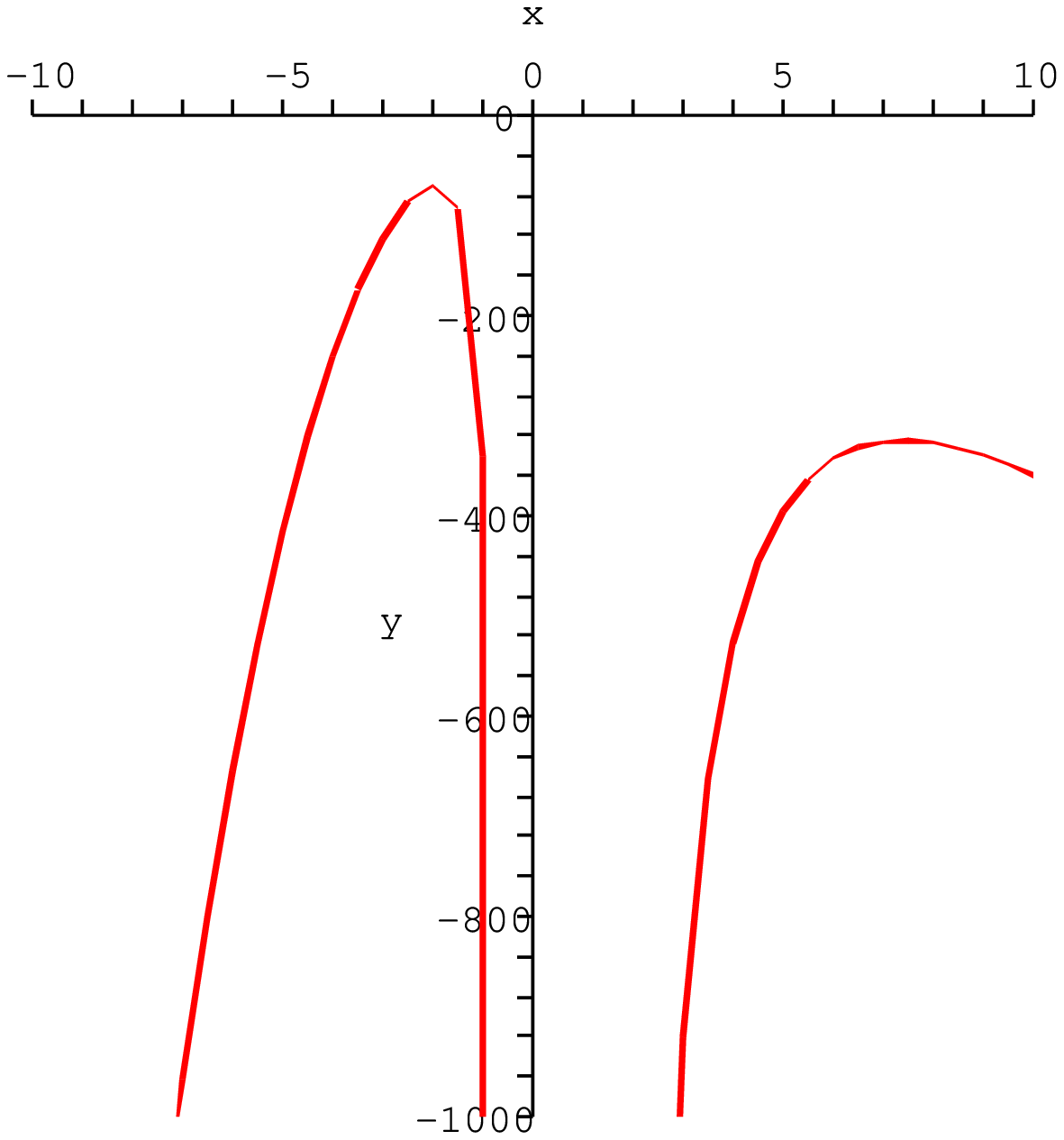}}}
\vskip 1cm
\caption{The curve $y^2=x^3+x$ and the  graph of the determinant of the three
tangent lines at $P_0+P,\ P_1+P, P_2+P$ computed in Maple.} \label{fig:counter}
\end{figure}

{\bf Example.} It is well known that any degree $3$ meromorphic
function on an elliptic curve $\C$ is realized in $\cP^2$ by the
standard equation $y^2=P_3(x)$ where $P_3(x)$ is a cubic
polynomial can be represented as the composition of the group
shift of the whole $\C$ by some fixed point on it  with the
projection from some point on $\cP^2$, see e.g. ??. The critical
points  on $\C\subset \cP^2$  for the projection from some point
$pt\in\cP^2$ are the points where the pencil of lines through $pt$
is tangent to $\C$. Thus in order to prove that a given triple
$(P_0,P_1,P_2)$ of real points on a given real elliptic curve
$\C\subset \cP^2$ can not serve as the set of critical points of a
real degree $3$ function we have to show that for any  choice of a
fourth real point $P\in \C$ the (real) tangent lines to $\C$ at
the points $P_0+P, P_1+P, P_2+P$ never meet at the same (real)
point on $\cP^2$. Recall that (up to a sign change) the addition
of two point $A$ and $B$ on a real elliptic $\C\subset \cP^2$ can
be interpreted as the third intersection point of the line
$\overline{AB}$ with  $\C$.

Consider the real elliptic curve $(\C,\sigma)$ given by the
equation $y^2=x^3+x$  with $3$  real points $P_0,P_1,P_2$, where $P_0$ is "almost"
infinite point $(100.35, 1005.3)$, $P_1$ is the point slightly
above the origin $(0.05,0.224)$, and, finally, $P_2$ is the
inflection point $(0.4,-0.67)$.

Then we claim that for any point $P$ on the real part $\C_\sigma$
there is no real function whose inflection points coincide with
$P_0+P, P_1+P, P_2+P$. Indeed as we explained above the existence of such a function
would mean that tangents to $\C_\sigma$ would intersects at the
same point of $\rP^2$, or the corresponding dual points in
$(\rP^2)^*$ are collinear. However, computing the determinant of
these three dual points we see that the obtained function of $P$
is non-vanishing, see Fig.~\ref{fig:counter}.

\section {Remarks and problems}
\label{sc:rmk}

\noindent I.  Analogously to the  total reality property for
rational curves one can ask a similar question for projective
curves of any genus, namely

\begin{problem}
    Given a real algebraic curve $(\C,\sigma)$ with compact $\C$ and nonempty real part
    $\C_{\sigma}$ and a complex algebraic map $\Psi:\C\to\cP^n$ such that the
inverse images of all the flattening points of $\Psi(\C)$ lie on the real
part  $\C_{\sigma}\subset \C$ is it true that  $\Psi$  is a real algebraic
 up to a M\"obius  transformation of the image $\cP^n$?
  \label{conj:extSS}
    \end{problem}

\noindent
    (Recent private communication from Gabrielov claims that this is now
    proven for the plane rational quintics.)

%
%
%
%
%
%

\medskip
 \noindent
 II. In the recent \cite{EGSV} the authors found another generalization of
 the conjecture on total reality in case of the usual rational functions.

 \begin{problem}  Extend the results of \cite{EGSV} to  the case of meromorphic
 functions  on curves of
 higher genera.
 \end{problem}

\end{document}